\documentclass[12pt]{article}
\usepackage{amssymb,amsfonts,amsmath}
\newcommand{\Z}{\mathbb Z}
\newcommand{\Q}{\mathbb Q}
\newcommand{\R}{\mathbb R}
\newcommand{\C}{\mathbb C}
\newcommand{\G}{\textbf{\textit{G}}}
\newcommand{\K}{\textbf{\textit{K}}}

\begin{document}
\thispagestyle{empty}
\hfill
\footnotetext{
\footnotesize
{\bf 2000 Mathematics Subject Classification}:~03B30,~12D99,~12L12,~14P05,~15A06,~26C10.
{\bf Key words and phrases:} system of polynomial (linear) equations, solution with minimal \mbox{$l\sb{\infty}$ norm.}}
\par
\noindent
\centerline{{\large Bounds of some real (complex) solution of a finite system}}
\par
\noindent
\centerline{{\large of polynomial equations with rational coefficients}}
\par
\noindent
\centerline{{\large Apoloniusz Tyszka}}
\par
\noindent
{\bf Abstract.} We discuss two conjectures. {\bf (I)} For each $x_1,\ldots,x_n \in \R$~($\C$)
there exist $y_1,\ldots,y_n \in \R$~($\C$) such that
\par
\noindent
\centerline{$\forall i \in \{1,\ldots,n\} ~|y_i| \leq 2^{\textstyle 2^{n-2}}$}
\par
\noindent
\centerline{$\forall i \in \{1,\ldots,n\} ~(x_i=1 \Rightarrow y_i=1)$}
\par
\noindent
\centerline{$\forall i,j,k \in \{1,\ldots,n\} ~(x_i+x_j=x_k \Rightarrow y_i+y_j=y_k)$}
\par
\noindent
\centerline{$\forall i,j,k \in \{1,\ldots,n\} ~(x_i \cdot x_j=x_k \Rightarrow y_i \cdot y_j=y_k)$}
\par
\noindent
{\bf (II)} Let $\G$ be an additive subgroup of $\C$. Then for each $x_1,\ldots,x_n \in \G$ there exist
$y_1,\ldots,y_n \in \G \cap \Q$ such that
\par
\noindent
\centerline{$\forall i \in \{1,\ldots,n\} ~|y_i| \leq 2^{n-1}$}
\par
\noindent
\centerline{$\forall i \in \{1,\ldots,n\} ~(x_i=1 \Rightarrow y_i=1)$}
\par
\noindent
\centerline{$\forall i,j,k \in \{1,\ldots,n\} ~(x_i+x_j=x_k \Rightarrow y_i+y_j=y_k)$}
\vskip 0.4truecm
\par
For a positive integer $n$ we define the set of equations $E_n$ by
\par
\noindent
\centerline{$E_n=\{x_i=1:~1 \leq i \leq n\}~\cup$}
\par
\noindent
\centerline{$\{x_i+x_j=x_k:~1 \leq i \leq j \leq n,~1 \leq k \leq n\}
\cup
\{x_i \cdot x_j=x_k:~1 \leq i \leq j \leq n,~1 \leq k \leq n\}$}
\par
\noindent
Since there is a finite number of non-empty subsets of~$E_n$, we get:
\begin{description}
\item{{\bf (1)}} There is a function $\chi:\{1,2,3,\ldots\} \to \{1,2,3,\ldots\}$ with the property:
for each positive integer~$n$, if a non-empty subset of~$E_n$ forms a system of equations that
is consistent over~$\Z$, then this system has a solution being
a sequence of integers whose absolute values are not greater than~$\chi(n)$.
\item{{\bf (2)}} There is a function $\gamma:\{1,2,3,\ldots\} \to \{1,2,3,\ldots\}$ with the property:
for each positive integer~$n$, if a non-empty subset of~$E_n$ forms a system of equations that
is consistent over~$\R$, then this system has a solution being
a sequence of real numbers whose absolute values are not greater than~$\gamma(n)$.
\item{{\bf (3)}} There is a function $\theta:\{1,2,3,\ldots\} \to \{1,2,3,\ldots\}$ with the property:
for each positive integer~$n$, if a non-empty subset of~$E_n$ forms a system of equations that
is consistent over~$\C$, then this system has a solution being
a sequence of complex numbers whose absolute values are not greater than~$\theta(n)$.
\end{description}
\vskip 0.2truecm
\par
\noindent
{\bf Remark.} Let us consider the problem of finding a complex
solution of a polynomial system with $n$~variables and integer coefficients.
This problem reduces to the problem of finding a real solution of some polynomial
system with $2n$~variables and integer coefficients. Therefore, the problem
of consistency over~$\C$ of a polynomial system with $n$~variables and integer coefficients
reduces to the problem of consistency over~$\R$ of some polynomial system with $2n$~variables
and integer coefficients.
\vskip 0.2truecm
\par
Let us note three facts:
\par
\noindent
\begin{description}
\item{{\bf (4)}}
There is a finite number of non-empty subsets of $E_n$.
\item{{\bf (5)}}
There is an algorithm for quantifier elimination for $(\R,+,\cdot,0,1,=,\leq)$
(A.~Tarski and A.~Seidenberg, and later G.~E.~Collins with his cylindrical algebraic decomposition algorithm, see \cite{Basu}).
In particular, there is an algorithm checking the consistency over~$\R$ of each finite
system of polynomial equations with integer coefficients.
\item{{\bf (6)}}
Applying the cylindrical algebraic decomposition algorithm,
for each consistent finite system of polynomial equations with $n$~variables and integer coefficients,
we can determine $a>0$ such that $[-a,a]^n$ contains a solution.
\end{description}
By the Remark and facts {\bf (4)}, {\bf (5)}, {\bf (6)}, we can find computable~$\gamma$ and~$\theta$.
\vskip 0.2truecm
\par
There are known direct estimates which enable us to find computable $\gamma$.
Let $V \subseteq \R^n$ be a real algebraic variety given by the system
of equations $f_1=\ldots=f_m=0$, where $f_i \in {\Q}[x_1,\ldots,x_n]$~$(i=1,\ldots,m)$.
We denote by $L$ the maximum of the bit-sizes of the coefficients
of the system and set $d=\sum_{i=1}^m\limits {\rm deg}(f_i)$, $r={n+2d \choose n}$.
We recall ([1,~p.~245]) that the bit-size of a non-zero integer is the number of
bits in its binary representation. More precisely, the bit-size of $k \in \Z \setminus \{0\}$
is $\tau$ if and only if $2^{\tau-1} \leq |k| < 2^{\tau}$.
The bit-size of a rational number is the sum of the bit-sizes of its
numerator and denominator.
N.~N.~Vorobjov Jr. proved~(\cite{Vorobjov}) that there exists $(x_1,\ldots,x_n) \in V$
such that $|x_i|<2^{H(r,L)}$~$(i=1,\ldots,n)$, where $H$ is some polynomial not
depending on the initial system; for a simplified proof see [7,~Lemma~9,~p.~56].
A~more general bound follows from [1,~Theorem~13.15,~p.~476].
\vskip 0.2truecm
\par
Let $M$ be the maximum of the absolute values of the coefficients of
the polynomials $f_1(x_1,\ldots,x_n),\ldots,f_m(x_1,\ldots,x_n) \in {\Z}[x_1,\ldots,x_n]$.
Assume that the system
\vskip 0.1truecm
\par
\noindent
\centerline{$f_1(x_1,\ldots,x_n)= \ldots =f_m(x_1,\ldots,x_n)=0$}
\vskip 0.1truecm
\par
\noindent
is consistent over $\R$ ($\C$).
Let $d(i,j)$ be the degree of variable~$x_i$ in the polynomial $f_j(x_1,\ldots,x_n)$.
Assume that $d_i={\rm max}\{d(i,j):~1 \leq j \leq m\} \geq 1$ for each $i \in \{1,\ldots,n\}$.
Let ${\cal T}$ denote the family of all polynomials $W(x_1,\ldots,x_n) \in {\Z}[x_1,\ldots,x_n]$
for which all coefficients belong to the interval $[-M,M]$ and each variable $x_i$ has degree less than
or equal to $d_i$. Then, $\{x_1,\ldots,x_n\} \subseteq {\cal T}$ and
\vskip 0.1truecm
\par
\noindent
\centerline{${\rm card} ~{\cal T}=(2M+1)^{(d_1+1) \cdot \ldots \cdot (d_n+1)}$}
\vskip 0.2truecm
\par
To each polynomial that belongs to ${\cal T} \setminus \{x_1,\ldots,x_n\}$ we assign
a new variable $x_i$ with $i \in \{n+1,\ldots,(2M+1)^{(d_1+1) \cdot \ldots \cdot (d_n+1)}\}$.
Let ${\cal K}$ denote the family of all equations of the form
\vskip 0.1truecm
\par
\noindent
\centerline{$x_i=1$, $x_i+x_j=x_k$, $x_i \cdot x_j=x_k$ ($i,j,k \in \{1, \ldots, (2M+1)^{(d_1+1) \cdot \ldots \cdot (d_n+1)}\}$)}
\vskip 0.1truecm
\par
\noindent
which are polynomial identities in ${\Z}[x_1,\ldots,x_n]$.
Let $f_j(x_1,\ldots,x_n)=x_{q(j)}$,
where $j \in \{1, \ldots, m\}$ and $q(j) \in \{1,\ldots,(2M+1)^{(d_1+1) \cdot \ldots \cdot (d_n+1)}\}$.
The system
\vskip 0.1truecm
\par
\noindent
\centerline{$f_1(x_1,\ldots,x_n)= \ldots =f_m(x_1,\ldots,x_n)=0$}
\vskip 0.1truecm
\par
\noindent
can be equivalently write down as
\vskip 0.1truecm
\par
\noindent
\centerline{${\cal K} \cup \{x_{q(1)}+x_{q(1)}=x_{q(1)},\ldots,x_{q(m)}+x_{q(m)}=x_{q(m)}\}$}
\vskip 0.2truecm
\par
We have proved that introducing additional variables we can equivalently write down the system
\vskip 0.1truecm
\par
\noindent
\centerline{$f_1(x_1,\ldots,x_n)= \ldots =f_m(x_1,\ldots,x_n)=0$}
\vskip 0.1truecm
\par
\noindent
as a system of equations of the form $x_i=1$, $x_i+x_j=x_k$, $x_i \cdot x_j=x_k$,
where $i,j,k \in \{1,\ldots,(2M+1)^{(d_1+1) \cdot \ldots \cdot (d_n+1)}\}$ and
the variables $x_{{\textstyle n+1}},\ldots,x_{{\textstyle (2M+1)^{(d_1+1) \cdot \ldots \cdot (d_n+1)}}}$
are new.
\vskip 0.2truecm
\par
\noindent
{\bf Conjecture 1.} Let $S$ be a consistent system of equations in real (complex)
numbers $x_1,x_2,\ldots,x_n$, where each equation in $S$ is one
of the following three forms: $x_i=1$ or $x_i+x_j=x_k$ or $x_i \cdot x_j=x_k$.
Then $S$ has a real (complex) solution $(x_1,\ldots,x_n)$ in which $|x_i| \leq 2^{\textstyle 2^{n-2}}$
for each $i$.
\vskip 0.2truecm
\par
\noindent
Conjecture 1 implies that the system
\vskip 0.1truecm
\par
\noindent
\centerline{$f_1(x_1,\ldots,x_n)= \ldots =f_m(x_1,\ldots,x_n)=0$}
\vskip 0.15truecm
\par
\noindent
has a real (complex) solution $(x_1,\ldots,x_n)$ in which
$|x_i| \leq 2^{2^{(2M+1)^{(d_1+1) \cdot \ldots \cdot (d_n+1)}-2}}$ for \hbox{each $i$.}
This upper bound is rough because does not depend on the number
of equations. We describe a better bound that depends on $m$.
Let
$$
L=\{(s_1,\ldots,s_n) \in {\Z}^n:~~(0 \leq s_1 \leq d_1) \wedge \ldots \wedge (0 \leq s_n \leq d_n) \wedge (1 \leq s_1+\ldots+s_n)\}
$$
$$
f_j(x_1,\ldots,x_n)=a_j+
\sum_{(s_1,\ldots,s_n) \in L} a_j(s_1,\ldots,s_n) \cdot x_1^{s_1} \ldots x_n^{s_n}
$$
where $j \in \{1,\ldots,m\}$, $a_j \in \Z$, $a_j(s_1,\ldots,s_n) \in \Z$. Let
\mbox{${\cal L}=\{x_1^{s_1} \ldots x_n^{s_n}:~(s_1,\ldots s_n) \in L\}$.}
Of course, $\{x_1,\ldots,x_n\} \subseteq {\cal L}$.
We define the lexicographic order $\preceq$ on ${\cal L}$.
We will define new variables $x_i$.
\begin{description}
\item{\emph{Step 1.}}
To each integer in $[-M,M]$ we assign a separate \mbox{variable~$x_i$}.\\
In this step we introduce $2M+1$ new variables.
\item{\emph{Step 2.}}
To each monomial in ${\cal L} \setminus \{x_1,\ldots,x_n\}$ we assign a separate variable $x_i$.\\
In this step we introduce $(d_1+1) \cdot \ldots \cdot (d_n+1)-1-n$ new variables.
\item{\emph{Step 3.}}
To each monomial $a_j(s_1,\ldots,s_n) \cdot x_1^{s_1} \ldots x_n^{s_n}$~~($j \in \{1,\ldots,m\}$, $(s_1,\ldots,s_n) \in L$)
we assign a separate variable $x_i$.\\
In this step we introduce $m \cdot ((d_1+1) \cdot \ldots \cdot (d_n+1)-1)$ new variables.
\item{\emph{Step 4.}}
To each polynomial $a_j+\sum_{(t_1,\ldots,t_n) \preceq (s_1,\ldots,s_n)}\limits~~a_j(t_1,\ldots,t_n) \cdot x_1^{t_1} \ldots x_n^{t_n}$\\
($j \in \{1,\ldots,m\}$, $(s_1,\ldots,s_n) \in L$) we assign a separate variable $x_i$.\\
In this step we introduce $m \cdot ((d_1+1) \cdot \ldots \cdot (d_n+1)-1)$ new variables.
\end{description}
\par
\noindent
The total number of new variables is equal to
\par
\noindent
\centerline{$p=2(M-m)-n+(2m+1) \cdot (d_1+1) \cdot \ldots \cdot (d_n+1)$}
\par
\noindent
Without lost of generality we can assume that we have introduced
the variables $x_{n+1},\ldots,x_{n+p}$.
Let ${\cal H}$ denote the family of all equations of the form
\vskip 0.2truecm
\par
\noindent
\centerline{$x_i=1$, $x_i+x_j=x_k$, $x_i \cdot x_j=x_k$~~($i,j,k \in \{1,\ldots,n+p\}$)}
\vskip 0.2truecm
\par
\noindent
which are polynomial identities in ${\Z}[x_1,\ldots,x_n]$.
Let $f_j(x_1,\ldots,x_n)=x_{q(j)}$,
where $j \in \{1, \ldots, m\}$ and $q(j) \in \{1,\ldots,n+p\}$.
The system
\vskip 0.2truecm
\par
\noindent
\centerline{$f_1(x_1,\ldots,x_n)= \ldots =f_m(x_1,\ldots,x_n)=0$}
\vskip 0.2truecm
\par
\noindent
can be equivalently write down as
\vskip 0.2truecm
\par
\noindent
\centerline{${\cal H} \cup \{x_{q(1)}+x_{q(1)}=x_{q(1)},\ldots,x_{q(m)}+x_{q(m)}=x_{q(m)}\}$}
\vskip 0.2truecm
\par
\noindent
Conjecture 1 implies that the system
\vskip 0.2truecm
\par
\noindent
\centerline{$f_1(x_1,\ldots,x_n)= \ldots =f_m(x_1,\ldots,x_n)=0$}
\vskip 0.2truecm
\par
\noindent
has a real (complex) solution $(x_1,\ldots,x_n)$ in which
$|x_i| \leq 2^{2^{2(M-m)-2+(2m+1) \cdot (d_1+1) \cdot \ldots \cdot (d_n+1)}}$ for each $i$.
\vskip 0.2truecm
\par
Concerning Conjecture 1, for $n=1$ estimation by $2^{\textstyle 2^{n-2}}$ can be replaced
by estimation by $1$. For $n>1$ estimation by $2^{\textstyle 2^{n-2}}$
is the best estimation. Indeed, let $n>1$ and
$\widetilde{x_1}=1$, $\widetilde{x_2}=2^{\textstyle 2^0}$,
$\widetilde{x_3}=2^{\textstyle 2^1}$, \ldots, $\widetilde{x_n}=2^{\textstyle 2^{n-2}}$.
In any ring $\K$ of characteristic $0$, from the system of all
equations belonging to $E_n$ and are satisfied under the substitution
$[x_1 \to \widetilde{x_1},~\ldots,~x_n \to \widetilde{x_n}]$,
it follows that $x_1=\widetilde{x_1}$, \ldots, $x_n=\widetilde{x_n}$.
\par
\vskip 0.2truecm
If a system $S \subseteq E_1$ is consistent over $\C$,
then $S$ has a solution $\widehat{x_1} \in \{0,1\}$.
This proves Conjecture 1 for $n=1$.
If a system $S \subseteq E_2$ is consistent over $\C$,
then $S$ has a solution $(\widehat{x_1},\widehat{x_2}) \in \{(0,0),(0,1),(1,0),(\frac{1}{2},1),(1,\frac{1}{2}),(1,1),(1,2),(2,1)\}$.
This proves Conjecture 1 for $n=2$. Let
\begin{eqnarray*}
{\cal W}=\bigl\{
\left\{1\right\},
\left\{0\right\},
\left\{1,0\right\},
\left\{1,2\right\},
\left\{1,\frac{1}{2}\right\},
\left\{1,2,\frac{1}{2}\right\},
\left\{1,0,2\right\},
\left\{1,0,\frac{1}{2}\right\},\\
\left\{1,0,-1\right\},
\left\{1,2,-1\right\},
\left\{1,2,3\right\},
\left\{1,2,4\right\},
\left\{1,\frac{1}{2},-\frac{1}{2}\right\},
\left\{1,\frac{1}{2},\frac{1}{4}\right\},
\left\{1,\frac{1}{2},\frac{3}{2}\right\},\\
\left\{1,-1,-2\right\},
\left\{1,\frac{1}{3},\frac{2}{3}\right\},
\left\{1,2,\sqrt{2}\right\},
\left\{1,\frac{1}{2},\frac{1}{\sqrt{2}}\right\},
\left\{1,\sqrt{2},\frac{1}{\sqrt{2}}\right\},\\
\left\{1,\frac{\sqrt{5}-1}{2},\frac{\sqrt{5}+1}{2}\right\},
\left\{1,\frac{\sqrt{5}+1}{2},\frac{\sqrt{5}+3}{2}\right\},
\left\{1,\frac{-\sqrt{5}-1}{2},\frac{\sqrt{5}+3}{2}\right\}
\bigr\}
\end{eqnarray*}
If a system $S \subseteq E_3$ is consistent over $\R$,
then $S$ has a real solution $(\widehat{x_1},\widehat{x_2},\widehat{x_3})$ with $\{\widehat{x_1}\} \cup \{\widehat{x_2}\} \cup \{\widehat{x_3}\} \in {\cal W}$.
This proves Conjecture 1 for $\R$ and $n=3$.
If a system $S \subseteq E_3$ is consistent over $\C$,
then $S$ has a solution $(\widehat{x_1},\widehat{x_2},\widehat{x_3})$ with $\{\widehat{x_1}\} \cup \{\widehat{x_2}\} \cup \{\widehat{x_3}\} \in
{\cal W} \cup \left\{\left\{1,\frac{-1+\sqrt{-3}}{2},\frac{1+\sqrt{-3}}{2}\right\},
\left\{1,\frac{1-\sqrt{-3}}{2},\frac{1+\sqrt{-3}}{2}\right\}\right\}$.
This proves Conjecture 1 for $\C$ and $n=3$.
\vskip 0.2truecm
\par
Now we demonstrate the use of the {\sl Mathematica}
software for checking \hbox{Conjecture 1} for $n=3$.
Without lost of generality, we can adopt the following
assumptions which allow for reducing the number of studied
systems $S \subseteq E_3$.
\begin{description}
\item{{\bf (7)}}
The equation $x_1=1$ belongs to system $S$,
as when all equations $x_1=1$, $x_2=1$, $x_3=1$ do not
belong to system $S$, then system $S$ has the solution $(0,0,0)$.
\item{{\bf (8)}}
Equations $x_1+x_1=x_1$, $x_1+x_2=x_2$, $x_1+x_3=x_3$
do not belong to system $S$, as each of these
equations is contradictory when $x_1=1$.
\item{{\bf (9)}}
We only consider these systems $S$ where each real
solution $(1, \widehat{x_2}, \widehat{x_3})$ has pairwise different
$1$, $\widehat{x_2}$, $\widehat{x_3}$, as Conjecture 1 is proven for $n=2$.
Therefore, all equations $x_1 \cdot x_i=x_j$ ($i \neq j$)
do not belong to system $S$.
\item{{\bf (10)}}
Instead of each equation $x_1+x_1=x_i$ ($i=2,3$)
we consider the equation $x_i=2$.
\item{{\bf (11)}}
Instead of each equation $x_i+x_i=x_1$ ($i=2,3$)
we consider the equation $x_i=\frac{1}{2}$.
\item{{\bf (12)}}
Instead of each equation $x_i+x_j=x_j$
($2 \leq i \leq j \leq 3$) we consider the equation $x_i=0$.
\item{{\bf (13)}}
Instead of each equation $x_i+x_j=x_i$
($1 \leq i \leq j \leq 3$, $2 \leq j$) we consider
the equation $x_j=0$.
\item{{\bf (14)}}
All equations $x_i \cdot x_j=x_i$
($1 \leq i \leq j \leq 3$) do not belong to system $S$,
because they are met when $j=1$, and when $j>1$ they are
equivalent to equation $x_i=0$.
\item{{\bf (15)}}
All equations $x_i \cdot x_j=x_j$
($1 \leq i \leq j \leq 3$) do not belong to system $S$,
because they are met when $i=1$, and when $i>1$ they are
equivalent to equation $x_j=0$.
\end{description}
After replacement of variables $x_1$, $x_2$, $x_3$
with $1$ and variables $x$, $y$, instead of set $E_3$
we receive 16 equations:
\begin{displaymath}
\begin{array}{cccccc}
x=2 & y=2 & x=\frac{1}{2} & y=\frac{1}{2} & x=0 & y=0 \\
x \cdot x=y & x \cdot x=1 & x+x=y & y \cdot y=x & y \cdot y=1 & y+y=x \\
& x \cdot y=1 & x+y=1 & x+1=y & y+1=x &
\end{array}
\end{displaymath}
\vskip 0.2truecm
\par
\noindent
If $n=3$, the following code in {\sl Mathematica} verifies
Conjecture 1 for $\R$ and $\C$:
\begin{quote}
\begin{verbatim}
Clear[x, y, i, j]
A := {x == 2, y == 2, x == 1/2, y == 1/2, x == 0, y == 0, x*x == y,
  x*x == 1, x + x == y, y*y == x, y*y == 1, y + y == x, x*y == 1,
  x + y == 1, x + 1 == y, y + 1 == x}
f[i_, j_] :=
 Reduce[Exists[{x, y},
    A[[i]] && A[[j]] && (Abs[x] > 4 || Abs[y] > 4)], Complexes] /;
  i < j
f[i_, j_] := {} /; i >= j
Union[Flatten[Table[f[i, j], {i, 1, 16}, {j, 1, 16}]]]
\end{verbatim}
\end{quote}
\par
\noindent
The output is {\tt \{False\}}.
\newpage
\par
\noindent
{\bf Theorem 1.} If $n \in \{1,2,3\}$, then Conjecture 1 holds true for each subring $\K \subseteq \C$.
\vskip 0.1truecm
\par
\noindent
{\it Proof.}
If a system $S \subseteq E_1$ is consistent over $\K$,
then $S$ has a solution $\widehat{x_1} \in \{0,1\}$.
If a system $S \subseteq E_2$ is consistent over $\K$ and $\frac{1}{2} \not\in \K$,
then $S$ has a solution $(\widehat{x_1},\widehat{x_2}) \in \{(0,0),(0,1),(1,0),(1,1),(1,2),(2,1)\}$.
If a system $S \subseteq E_2$ is consistent over $\K$ and $\frac{1}{2} \in \K$,
then $S$ has a solution $(\widehat{x_1},\widehat{x_2}) \in \{(0,0),(0,1),(1,0),(\frac{1}{2},1),(1,\frac{1}{2}),(1,1),(1,2),(2,1)\}$.
For reducing the number of studied systems $S \subseteq E_3$,
we may assume that the equation $x_1=1$ belongs to $S$,
as when all equations $x_1=1$, $x_2=1$, $x_3=1$ do not belong to~$S$,
then $S$ has the solution $(0,0,0) \in {\K}^3$. Let
\vskip 0.1truecm
\par
\noindent
\centerline{$A_2=\{\widehat{x_2} \in \C: {\rm ~there~exists~} \widehat{x_3} \in \C {\rm ~for~which~} (1,\widehat{x_2},\widehat{x_3}) {\rm ~solves~} S\}$}
\vskip 0.1truecm
\par
\noindent
\centerline{$A_3=\{\widehat{x_3} \in \C: {\rm ~there~exists~} \widehat{x_2} \in \C {\rm ~for~which~} (1,\widehat{x_2},\widehat{x_3}) {\rm ~solves~} S\}$}
\vskip 0.1truecm
\par
\noindent
We may assume that $A_2 \not\subseteq \{z \in \C: |z| \leq 4\}$
or $A_3 \not\subseteq \{z \in \C: |z| \leq 4\}$.
\vskip 0.1truecm
\par
\noindent
Case 1: $A_2 \not\subseteq \{z \in \C: |z| \leq 4\}$ and
$A_3 \subseteq \{z \in \C: |z| \leq 4 \}$.
If $(1,\widehat{x_2},\widehat{x_3}) \in {\K}^3$ solves~$S$,
then $(1,1,\widehat{x_3}) \in {\K}^3$ solves~$S$.
\vskip 0.1truecm
\par
\noindent
Case 2: $A_2 \subseteq \{z \in \C: |z| \leq 4\}$ and
$A_3 \not\subseteq \{z \in \C: |z| \leq 4 \}$.
If $(1,\widehat{x_2},\widehat{x_3}) \in {\K}^3$ solves~$S$,
then $(1,\widehat{x_2},1) \in {\K}^3$ solves~$S$.
\vskip 0.1truecm
\par
\noindent
Case 3: $A_2 \not\subseteq \{z \in \C: |z| \leq 4\}$ and
$A_3 \not\subseteq \{z \in \C: |z| \leq 4 \}$.
If $(1,\widehat{x_2},\widehat{x_3}) \in {\K}^3$ solves~$S$,
then $(1,0,1) \in {\K}^3$ solves~$S$ or $(1,1,0) \in {\K}^3$ solves~$S$ or $(1,1,1) \in {\K}^3$ solves~$S$.
\newline
\rightline{$\Box$}
\vskip 0.2truecm
\par
Conjecture 1 holds true for each $n \in \{1,2,3,4\}$ and each subring $\K \subseteq \C$.
It follows from the following Observation 1 which borrows the idea from the proof of Theorem 1.
\vskip 0.2truecm
\par
\noindent
{\bf Observation 1.} Let $n \in \{1,2,3,4\}$, and let $S \subseteq E_n$ be a system that is consistent over the subring $\K \subseteq \C$.
If $(x_1,\ldots,x_n) \in {\K}^n$ solves $S$, then $(\widehat{x_1},\ldots,\widehat{x_n})$ solves $S$,
where each $\widehat{x_i}$ is suitably chosen from $\{x_i,0,1,2,\frac{1}{2}\} \cap \{z \in \K:~|z| \leq 2^{\textstyle 2^{n-2}}\}$.
\vskip 0.2truecm
\par
Multiple execution of the algorithm described in items
{\bf (16)}-{\bf (20)} yields partial (as probabilistic)
resolution of Conjecture 1 for $\R$ and $n \geq 4$.
This algorithm resolves Conjecture 1 for some randomly chosen subsystem of $E_n$.
\newpage
\par
\noindent
{\bf (16)}~~From the set $E_n$ we remove the equations
\vskip 0.1truecm
\par
\noindent
\centerline{$
\begin{array}{ll}
x_i=1             & (1 \leq i \leq n)\\
x_1+x_i=x_i       & (1 \leq i \leq n)\\
x_1+x_1=x_i       & (i \neq 2)\\
x_i+x_i=x_1       & (i \neq 3)\\
x_i+x_j=x_i       & (1 \leq i \leq j \leq n, (i,j) \neq (4,4))\\
x_i+x_j=x_j       & (1 \leq i \leq j \leq n, (i,j) \neq (4,4))\\
x_i \cdot x_j=x_i & (1 \leq i \leq j \leq n)\\
x_i \cdot x_j=x_j & (1 \leq i \leq j \leq n)\\
x_1 \cdot x_i=x_j & (i,j \in \{1,2,\ldots,n\})
\end{array}$}
\vskip 0.1truecm
\par
\noindent
and in other equations we replace $x_1$ by number $1$. We obtain a non-empty set $H_n$.
\vskip 0.2truecm
\par
\noindent
{\bf (17)}~~We introduce a random linear order on $H_n$,
but with a reservation that the first equation is to be
among equations involving number $1$.
\vskip 0.2truecm
\par
\noindent
{\bf (18)}~~We define by induction a finite sequence $(s_1,\ldots,s_m)$
of equations belonging \hbox{to $H_n$.} As $s_1$ we put the first equation in $H_n$.
After this, we remove from $H_n$ all equations having the left side identical
to the left side of equation $s_1$; this step may be omitted.
When the sequence $(s_1,\ldots,s_i)$ is defined,
and there exists \hbox{$h \in H_n \setminus \{s_1,\ldots,s_i\}$}
for which the system $\{s_1,\ldots,s_i, h\}$ has a real solution
$(x_2,\ldots,x_n)$ with pairwise different
$1,x_2,\ldots,x_n$, then as $s_{i+1}$ we put the smallest such $h$.
After this, we remove from $H_n$ all equations having the left side identical
to the left side of equation $s_{i+1}$; this step may be omitted.
If such $h$ does not exist, then $m=i$ and the construction of the sequence
$(s_1,\ldots,s_m)$ is finished. 
\vskip 0.2truecm
\par
The condition "with pairwise different $1,x_2,\ldots,x_n$" may be removed from \hbox{item {\bf (18)},}
but this will increase the average number of executions of item {\bf (18)}.
\vskip 0.2truecm
\par
\noindent
{\bf (19)}~~If any of the systems
$\{s_1,\ldots,s_m,x_i=1\}$~$(i=2,\ldots,n)$,
$\{s_1,\ldots,s_m,x_i=x_j\}$~$(2 \leq i<j \leq n)$ has a real
solution, then we return to item {\bf (17)}.
\vskip 0.2truecm
\par
\noindent
{\bf (20)}~~We resolve Conjecture 1 for $\R$ for the latest system $\{s_1,\ldots,s_m\}$.
\vskip 0.2truecm
\par
The above algorithm resolves Conjecture 1 only for these subsystems of $H_n$
for which each real solution $(x_2,\ldots,x_n)$
satisfies ${\rm card}(\{1,x_2,\ldots,x_n\})=n$.
It is sufficient if Conjecture 1 was previously resolved for $n-1$ real variables.
\newpage
In the computer execution of items {\bf (18)}-{\bf (20)}
one may use {\sl Mathematica} software and the
{\sl CylindricalDecomposition} or {\sl Reduce} procedure.
The algorithm for \hbox{Conjecture 1} for $\C$ is analogical,
but we only apply the {\sl Reduce} procedure.
Neither for $\R$ nor for $\C$ can we apply the {\sl Resolve}
procedure in Mathematica 6.0.1, as it yields wrong results, as presented below:
\begin{quote}
\begin{verbatim}
Resolve[Exists[{x}, x == 0 && x*x == 1], Reals]
True
Resolve[Exists[{x}, x == 0 && x*x == 1], Complexes]
True
\end{verbatim}
\end{quote}
\vskip 0.2truecm
\par
The number $2+273^2$ is prime.
\vskip 0.2truecm
\par
\noindent
{\bf Theorem 2.} If $k \in \Z \cap [273,\infty)$ and $2+k^2$ is prime,
then Conjecture 1 fails for $n=6$ and the ring ${\Z}[\frac{1}{2+k^2}]=\{\frac{x}{(2+k^2)^m}:~x \in \Z,m \in \Z \cap [0,\infty)\}$.
\vskip 0.2truecm
\par
\noindent
{\it Proof.} $(1,2,k,k^2,2+k^2,\frac{1}{2+k^2})$ solves the system
\begin{displaymath}
\left\{
\begin{array}{rcl}
x_1 &=& 1 \\
x_1+x_1 &=& x_2 \\
x_3 \cdot x_3 &=& x_4 \\
x_2+x_4 &=& x_5 \\
x_5 \cdot x_6 &=& x_1
\end{array}
\right.
\end{displaymath}
Assume that $(x_1,x_2,x_3,x_4,x_5,x_6) \in ({\Z}[\frac{1}{2+k^2}])^6$
solves the system. Let $x_5=\frac{a}{(2+k^2)^p}$, $x_6=\frac{b}{(2+k^2)^q}$,
$a,b \in \Z$, $p,q \in \Z \cap [0,\infty)$. Since $2+k^2$ is prime and
$1=|x_1|=|x_5 \cdot x_6|=\frac{|a| \cdot |b|}{(2+k^2)^{p+q}}$, we conclude that
$|a|=(2+k^2)^{\widetilde{p}}$ for some $\widetilde{p} \in \Z \cap [0,\infty)$.
Hence $|x_5|=(2+k^2)^{\widetilde{p}-p}$. On the other hand,
$|x_5|=|x_2+x_4|=|x_1+x_1+x_3 \cdot x_3|=|1+1+x_3^2| \geq 2$. Therefore, $\widetilde{p}-p \geq 1$.
Consequently, $|x_5|=(2+k^2)^{\widetilde{p}-p} \geq 2+k^2>2^{\textstyle 2^{6-2}}$.
\newline
\rightline{$\Box$}
\par
Lemma 1 is a special case of the result presented in [16,~p.~3].
\par
\noindent
{\bf Lemma 1.} For each non-zero integer $x$ there exist integers $a$, $b$ such that $ax=(2b-1)(3b-1)$.
\par
\noindent
{\it Proof.} Write $x$ as $(2y-1) \cdot 2^m$, where $y \in \Z$ and $m \in \Z \cap [0,\infty)$.
Obviously, \hbox{$\frac{\textstyle 2^{2m+1}+1}{\textstyle 3} \in \Z$.}
By Chinese Remainder Theorem we can find an integer $b$ such that
\hbox{$b \equiv y {\rm ~(mod~} 2y-1)$}
and \hbox{$b \equiv \frac{\textstyle 2^{2m+1}+1}{\textstyle 3} {\rm ~(mod~} 2^m)$.}
Thus, \hbox{$\frac{\textstyle 2b-1}{\textstyle 2y-1} \in \Z$} and
\hbox{$\frac{\textstyle 3b-1}{\textstyle 2^m} \in \Z$.}
Hence $\frac{\textstyle (2b-1)(3b-1)}{\textstyle x}=\frac{\textstyle 2b-1}{\textstyle 2y-1} \cdot \frac{\textstyle 3b-1}{\textstyle 2^m} \in \Z$.
\newline
\rightline{$\Box$}
\par
\noindent
{\bf Theorem 3.} If a prime number $p$ is greater than $2^{256}$, then
Conjecture 1 fails for $n=10$ and the ring ${\Z}\left[\frac{1}{p}\right]$.
\vskip 0.2truecm
\par
\noindent
{\it Proof.} Let us consider the system
\begin{displaymath}
\left\{
\begin{array}{rcl}
x_1 &=& 1\\
x_2 \cdot x_3 &=& x_1 \\
x_3+x_4 &=& x_2 \\
x_4 \cdot x_5 &=& x_6 \\
x_7+x_7 &=& x_8 \\
x_1+x_9 &=& x_8 \\
x_7+ x_9 &=& x_{10} \\
x_9 \cdot x_{10} &=& x_6
\end{array}
\right.
\end{displaymath}
By Lemma 1 there exist integers $u$, $s$ such that $(p^2-1) \cdot u=(2s-1)(3s-1)$.
Hence
\vskip 0.2truecm
\par
\noindent
\centerline{$(1,~p,~\frac{1}{p},~p-\frac{1}{p},~p \cdot u,~(p^2-1) \cdot u,~s,~2s,~2s-1,~3s-1) \in \left({\Z}\left[\frac{1}{p}\right]\right)^{10}$}
\vskip 0.2truecm
\par
\noindent
solves the system. If $(x_1,x_2,x_3,x_4,x_5,x_6,x_7,x_8,x_9,x_{10}) \in \left({\Z}\left[\frac{1}{p}\right]\right)^{10}$
solves the system, then $(x_2-x_3) \cdot x_5=(2x_7-1)(3x_7-1)$.
Since $2x_7-1 \neq 0$ and $3x_7-1 \neq 0$, we get $x_2 \neq x_3$.
Since $x_2 \cdot x_3=1$, we get:~$|x_2|=p^n$ for some $n \in \Z \cap [1,\infty)$ or $|x_3|=p^n$
for some $n \in \Z \cap [1,\infty)$. Therefore, $|x_2| \geq p>2^{\textstyle 2^{10-2}}$ or
$|x_3| \geq p>2^{\textstyle 2^{10-2}}$.
\newline
\rightline{$\Box$}
\vskip 0.2truecm
\par
The number $-2^{32}-2^{16}-1$ is square-free, because
$-3 \cdot 7 \cdot 13 \cdot 97 \cdot 241 \cdot 673$
is the factorization of $-2^{32}-2^{16}-1$ into prime numbers.
\vskip 0.2truecm
\par
\noindent
{\bf Theorem 4.} Conjecture 1 fails for $n=6$ and the ring
${\Z}[\sqrt{-2^{32}-2^{16}-1}]=\{x+y \cdot \sqrt{-2^{32}-2^{16}-1}:~x,y \in \Z\}$.
\vskip 0.2truecm
\par
\noindent
{\it Proof.} $(1,2^{16}+1,-2^{16},-2^{32}-2^{16},\sqrt{-2^{32}-2^{16}-1},-2^{32}-2^{16}-1)$ solves the system
\begin{displaymath}
\left\{
\begin{array}{rcl}
x_1 &=& 1 \\
x_2+x_3 &=& x_1 \\
x_2 \cdot x_3 &=& x_4 \\
x_5 \cdot x_5 &=& x_6 \\
x_1+x_6 &=& x_4
\end{array}
\right.
\end{displaymath}
which has no integer solutions. For each $z \in {\Z}[\sqrt{-2^{32}-2^{16}-1}]$,
if $|z| \leq 2^{\textstyle 2^{6-2}}$ then $z \in \Z$.
\newline
\rightline{$\Box$}
\vskip 0.2truecm
\par
\noindent
{\bf Observation 2.} If $q,a,b,c,d \in \Z$, $b \neq 0$ or $d \neq 0$, $q \geq 2$, $q$ is square-free,
and $(a+b\sqrt{q}) \cdot (c+d\sqrt{q})=1$, then
$$
(a \geq 1 \wedge b \geq 1) \vee (a \leq -1 \wedge b \leq -1) \vee (c \geq 1 \wedge d \geq 1) \vee (c \leq -1 \wedge d \leq -1)
$$
\par
The number $4 \cdot 13^4-1$ is square-free, because $3 \cdot 113 \cdot 337$ is
the factorization of $4 \cdot 13^4-1$ into prime numbers.
\vskip 0.2truecm
\par
\noindent
{\bf Theorem 5.} If $p \in \Z \cap [13,\infty)$ and $4p^4-1$ is square-free,
then Conjecture 1 fails for $n=5$ and the ring ${\Z}[\sqrt{4p^4-1}]=\{x+y \cdot \sqrt{4p^4-1}:~x,y \in \Z\}$.
\vskip 0.2truecm
\par
\noindent
{\it Proof.} $(1,2p^2+\sqrt{4p^4-1},2p^2-\sqrt{4p^4-1},4p^2,2p)$ solves the system
\begin{displaymath}
\left\{
\begin{array}{rcl}
x_1 &=& 1 \\
x_2 \cdot x_3 &=& x_1 \\
x_2+x_3 &=& x_4 \\
x_5 \cdot x_5 &=& x_4
\end{array}
\right.
\end{displaymath}
Assume that $(x_1,x_2,x_3,x_4,x_5) \in ({\Z}[\sqrt{4p^4-1}])^5$ solves the system.
Let $x_2=a+b\sqrt{4p^4-1}$, $x_3=c+d\sqrt{4p^4-1}$, $a,b,c,d \in \Z$. Since
$$
\neg (\exists x_2 \in \Z ~\exists x_3 \in \Z ~\exists x_5 \in {\Z}[\sqrt{4p^4-1}] ~(x_2 \cdot x_3=1 \wedge x_2+x_3=x_5^2)),
$$
$b \neq 0$ or $d \neq 0$. Since $x_2\cdot x_3=1$, Observation 2 implies that
$|x_2| \geq 1+\sqrt{4p^4-1}>2^{\textstyle 2^{5-2}}$ or $|x_3| \geq 1+\sqrt{4p^4-1}>2^{\textstyle 2^{5-2}}$.
\newline
\rightline{$\Box$}
\vskip 0.2truecm
\par
\noindent
{\bf Theorem 6.} Let $f(x,y) \in {\Q}[x,y]$ and the equation $f(x,y)=0$ defines an
irreducible algebraic curve of genus greater than $1$.
Let some $r \in \R$ satisfies
$$
(\ast)~~(-\infty,r) \subseteq \{x \in \R:~\exists y \in \R ~f(x,y)=0\}~ \vee ~(r,\infty) \subseteq \{x \in \R:~\exists y \in \R ~f(x,y)=0\}
$$
and let $\K$ denote the function field over $\Q$ defined by $f(x,y)=0$.
Then \hbox{Conjecture 1} fails for some subfield of $\R$ that is isomorphic to $\K$.
\vskip 0.2truecm
\par
\noindent
{\it Proof.} By Faltings' finiteness theorem (\cite{Faltings},~cf.~[11,~p.~12]) the set
$$
\{u \in \K:~ \exists v \in \K ~f(u,v)=0\}
$$
is finite. Let ${\rm card~}\{u \in \K:~ \exists v \in \K ~f(u,v)=0\}=n \geq 1$,
and let ${\cal U}$ denote the following system of equations
\begin{displaymath}
\left\{
\begin{array}{rcl}
f(x_i,y_i) &=& 0 ~~(1 \leq i \leq n) \\
x_i+t_{i,j} &=& x_j ~~~(1 \leq i<j \leq n) \\
t_{i,j} \cdot s_{i,j} &=& 1 ~~(1 \leq i<j \leq n) \\
x_{n+1} &=& \sum_{i=1}^n\limits x_i^2
\end{array}
\right.
\end{displaymath}
For some integer $m>n$ there exists a set ${\cal G}$ of $m$ variables such that
$$
\{x_1,\ldots,x_n\,x_{n+1},y_1,\ldots,y_n\} \cup \{t_{i,j},s_{i,j}:~1 \leq i<j \leq n\} \subseteq {\cal G}
$$
and the system ${\cal U}$ can be equivalently write down as a system ${\cal V}$ which contains only equations
of the form $X=1$, $X+Y=Z$, $X \cdot Y=Z$, where $X,Y,Z \in {\cal G}$.
By ($\ast$), we find $\widetilde{x},\widetilde{y} \in \R$ such that $f(\widetilde{x},\widetilde{y})=0$,
$\widetilde{x}$ is transcendental over $\Q$, and $|\widetilde{x}|>2^{\textstyle 2^{m-3}}$.
If $(\widehat{x_1}, \ldots, \widehat{x_m}) \in ({\Q}(\widetilde{x},\widetilde{y}))^m$ solves ${\cal V}$, then
$$
\widehat{x_{n+1}} = \sum_{i=1}^n \widehat{x_i}^2 \geq \widetilde{x}^2 > (2^{\textstyle 2^{m-3}})^2 = 2^{\textstyle 2^{m-2}}
$$
Obviously, $\K$ is isomorphic to ${\Q}(\widetilde{x},\widetilde{y})$.
\newline
\rightline{$\Box$}
\vskip 0.2truecm
\par
Conjecture 1 fails for some subfield of $\R$ and $n=7$.
We sketch the proof here.
We find $\alpha,\beta \in \R$ such that $\alpha^2 \cdot \beta \cdot (1-\alpha^2-\beta)=1$,
$\alpha$ is transcendental \hbox{over $\Q$,} and $|\alpha|>2^{\textstyle 2^{7-2}}$.
It is known (\cite{Sansone}) that the equation $x+y+z=xyz=1$ has no rational solution.
Applying this, we prove: if $(x_1,x_2,x_3,x_4,x_5,x_6,x_7) \in {\Q}(\alpha,\beta)^7$ solves the system
\begin{displaymath}
\left\{
\begin{array}{rcl}
x_1 &=& 1 \\
x_2 \cdot x_2 &=& x_3 \\
x_3 + x_4 &=& x_5 \\
x_5 + x_6 &=& x_1 \\
x_3 \cdot x_4 &=& x_7 \\
x_6 \cdot x_7 &=& x_1
\end{array}
\right.
\end{displaymath}
then $|x_2|=|\alpha|>2^{\textstyle 2^{7-2}}$.
\vskip 0.2truecm
\par
For each $a,b,c \in \R$~($\C$) we define $S(a,b,c)$ as
\vskip 0.2truecm
\par
\noindent
\centerline{$\{{\cal E} \in E_3:~{\cal E} {\rm ~is~satisfied~under~the~substitution~}[x_1 \to a,~x_2 \to b,~x_3 \to c]\}$}
\vskip 0.2truecm
\par
If $a,b,c \in \R$ and $\{a\} \cup \{b\} \cup \{c\} \in {\cal W}$, then
the system $S(a,b,c)$ is consistent over~$\R$, has a finite number
of real solutions, and each real solution of~$S(a,b,c)$ belongs
to $[-4,4]^3$. The family
$$
\{S(a,b,c):~ a,b,c \in \R ~\wedge~ \{a\} \cup \{b\} \cup \{c\} \in {\cal W}\}
$$
equals to the family of all systems $S \subseteq E_3$ which are
consistent over $\R$ and maximal with respect to inclusion.
\vskip 0.2truecm
\par
If $a,b,c \in \C$ and $\{a\} \cup \{b\} \cup \{c\} \in
{\cal W} \cup \left\{\left\{1,\frac{-1+\sqrt{-3}}{2},\frac{1+\sqrt{-3}}{2}\right\},
\left\{1,\frac{1-\sqrt{-3}}{2},\frac{1+\sqrt{-3}}{2}\right\}\right\}$, then
the system $S(a,b,c)$ is consistent over $\C$, has a finite number
of solutions, and each solution of $S(a,b,c)$ belongs to
$\{(z_1,z_2,z_3) \in {\C}^3:~|z_1| \leq 4 ~\wedge~ |z_2| \leq 4 ~\wedge~ |z_3| \leq 4\}$.
The family
\vskip 0.2truecm
\par
\noindent
\centerline{$\{S(a,b,c):~ a,b,c \in \C ~\wedge~ \{a\} \cup \{b\} \cup \{c\} \in$}
\vskip 0.2truecm
\par
\noindent
\centerline{${\cal W} \cup \left\{\left\{1,\frac{-1+\sqrt{-3}}{2},\frac{1+\sqrt{-3}}{2}\right\},
\left\{1,\frac{1-\sqrt{-3}}{2},\frac{1+\sqrt{-3}}{2}\right\}\right\}\}$}
\vskip 0.2truecm
\par
\noindent
equals to the family of all systems $S \subseteq E_3$ which
are consistent over $\C$ and maximal with respect to inclusion.
\vskip 0.2truecm
\par
Let us consider the following four conjectures, analogical conjectures seem to be true for $\R$.
\begin{description}
\item{{\bf (21a)}}
If a system $S \subseteq E_n$ is consistent over $\C$ and maximal with respect to inclusion,
then each solution of $S$ belongs to\\
$\{(x_1,\ldots,x_n) \in {\C}^n:~|x_1| \leq 2^{\textstyle 2^{n-2}}~\wedge~\ldots~\wedge~|x_n| \leq 2^{\textstyle 2^{n-2}}\}$.
\item{{\bf (21b)}}
If a system $S \subseteq E_n$ is consistent over $\C$ and maximal with respect to inclusion,
then $S$ has a finite number of solutions $(x_1,\ldots,x_n)$.
\item{{\bf (21c)}}
If the equation $x_1=1$ belongs to $S \subseteq E_n$ and $S$ has a finite number of complex solutions $(x_1,\ldots,x_n)$,
then each such solution belongs to\\
$\{(x_1,\ldots,x_n) \in {\C}^n:~|x_1| \leq 2^{\textstyle 2^{n-2}}~\wedge~\ldots~\wedge~|x_n| \leq 2^{\textstyle 2^{n-2}}\}$.
\item{{\bf (21d)}}
If a system $S \subseteq E_n$ has a finite number of complex solutions $(x_1,\ldots,x_n)$,
then each such solution belongs to\\
$\{(x_1,\ldots,x_n) \in {\C}^n:~|x_1| \leq 2^{\textstyle 2^{n-1}}~\wedge~\ldots~\wedge~|x_n| \leq 2^{\textstyle 2^{n-1}}\}$.
\end{description}
\par
Conjecture 21a strengthens Conjecture 1 for $\C$. The conjunction of Conjectures 21b and 21c implies Conjecture 21a.
\vskip 0.2truecm
\par
Concerning Conjecture 21d, for $n=1$ estimation by $2^{\textstyle 2^{n-1}}$ can be replaced
by estimation by $1$. For $n>1$ estimation by $2^{\textstyle 2^{n-1}}$ is the best estimation. Indeed, the system
\begin{displaymath}
\left\{
\begin{array}{rcl}
x_1+x_1 &=& x_2 \\
x_1 \cdot x_1 &=& x_2 \\
x_2 \cdot x_2 &=& x_3 \\
x_3 \cdot x_3 &=& x_4 \\
&...& \\
x_{n-1} \cdot x_{n-1} &=&x_n
\end{array}
\right.
\end{displaymath}
has precisely two complex solutions, $(0,\ldots,0)$, $(2,4,16,256,\ldots,2^{\textstyle 2^{n-2}},2^{\textstyle 2^{n-1}})$.
\vskip 0.2truecm
\par
The following code in {\sl MuPAD} yields a probabilistic confirmation of Conjectures 21b and 21c.
The value of $n$ is set, for example, to $5$. The number of iterations is set, for example, to $1000$.
\begin{quote}
\begin{verbatim}
SEED:=time():
p:=[v-1,x-1,y-1,z-1]:
var:=[1,v,x,y,z]:
for i from 1 to 5 do
for j from i to 5 do
for k from 1 to 5 do
p:=append(p,var[i]+var[j]-var[k]):
p:=append(p,var[i]*var[j]-var[k]):
end_for:
end_for:
end_for:
p:=listlib::removeDuplicates(p):
max_abs_value:=1:
for r from 1 to 1000 do
q:=combinat::permutations::random(p):
syst:=[t-v-x-y-z]:
w:=1:
repeat
if groebner::dimension(append(syst,q[w]))>-1
then syst:=append(syst,q[w]) end_if:
w:=w+1:
until (groebner::dimension(syst)=0 or w>nops(q)) end:
d:=groebner::dimension(syst):
if d>0 then print("Conjecture 21b is false") end_if:
if d=0 then
sol:=numeric::solve(syst):
for m from 1 to nops(sol) do
for n from 2 to 5 do
max_abs_value:=max(max_abs_value,abs(sol[m][n][2])):
end_for:
end_for:
end_if:
print(max_abs_value);
end_for:
\end{verbatim}
\end{quote}
\par
\noindent
If we replace
\begin{quote}
\begin{verbatim}
p:=[v-1,x-1,y-1,z-1]:    by    p:=[]:
var:=[1,v,x,y,z]:        by    var:=[u,v,x,y,z]:
max_abs_value:=1:        by    max_abs_value:=0:
syst:=[t-v-x-y-z]:       by    syst:=[t-u-v-x-y-z]:
for n from 2 to 5 do     by    for n from 2 to 6 do
\end{verbatim}
\end{quote}
then we get a code for a probabilistic confirmation of Conjecture 21d.
\vskip 0.2truecm
\par
We can formulate Conjecture 1 as follows:
for each $x_1,\ldots,x_n \in \R$~($\C$) there exist $y_1,\ldots,y_n \in \R$~($\C$) such that
\par
\noindent
\centerline{$\forall i \in \{1,\ldots,n\} ~|y_i| \leq 2^{\textstyle 2^{n-2}}$}
\par
\noindent
\centerline{$\forall i \in \{1,\ldots,n\} ~(x_i=1 \Rightarrow y_i=1)$}
\par
\noindent
\centerline{$\forall i,j,k \in \{1,\ldots,n\} ~(x_i+x_j=x_k \Rightarrow y_i+y_j=y_k)$}
\par
\noindent
\centerline{$\forall i,j,k \in \{1,\ldots,n\} ~(x_i \cdot x_j=x_k \Rightarrow y_i \cdot y_j=y_k)$}
\vskip 0.2truecm
\par
We say that $X \subseteq \R$~($\C$) has a property ${\cal B}$,
if for each $x_1,\ldots,x_n \in X$ there exist $y_1,\ldots,y_n \in X$
with the above four properties. We define:
$$
{\cal F}_{\R}=\{X \subseteq \R:~X{\rm ~has~property~} {\cal B}\}
$$
$$
{\cal F}_{\C}=\{X \subseteq \C:~X {\rm ~has~property~} {\cal B}\}
$$
If $X \subseteq [-\sqrt{2},~\sqrt{2}]$ then $X \in {\cal F}_{\R}$.
If $X \subseteq \{z \in \C:~ |z| \leq \sqrt{2}\}$ then $X \in {\cal F}_{\C}$.
\vskip 0.2truecm
\par
\noindent
{\bf Theorem 7.} The family ${\cal F}_{\R}$ (${\cal F}_{\C}$)
has a maximal element.
\vskip 0.2truecm
\par
\noindent
{\it Proof.} We prove: if ${\cal C} \subseteq {\cal F}_{\R}$ is a chain,
then $\bigcup_{X \in {\cal C}}\limits X \in {\cal F}_{\R}$.
Since ${\cal C}$ is a chain, for each $x_1,\ldots,x_n \in \bigcup_{X \in {\cal C}}\limits X$
there exists $X \in {\cal C}$ with $x_1,\ldots,x_n \in X$.
Since $X$ has property ${\cal B}$, we obtain
suitable $y_1,\ldots,y_n \in X \cap [-2^{\textstyle 2^{n-2}},~2^{\textstyle 2^{n-2}}]
\subseteq \left( \bigcup_{X \in {\cal C}}\limits X \right) \cap [-2^{\textstyle 2^{n-2}},~2^{\textstyle 2^{n-2}}]$.
We have proved that $\bigcup_{X \in {\cal C}}\limits X \in {\cal F}_{\R}$.
By Zorn's lemma, the family ${\cal F}_{\R}$ has a maximal element.
The proof for ${\cal F}_{\C}$ is analogical.
\newline
\rightline{$\Box$}
\par
It is hardly to decide whether Theorem 7 may help prove that $\R \in {\cal F}_{\R}$ and $\C \in {\cal F}_{\C}$.
\vskip 0.2truecm
\par
Conjecture 2 strengthens Conjecture 1.
\vskip 0.2truecm
\par
\noindent
{\bf Conjecture 2.} For each $x_1,\ldots,x_n \in \R$~($\C$) there exist $y_1,\ldots,y_n \in \R$~($\C$) such that
\par
\noindent
\centerline{$\forall i \in \{1,\ldots,n\} ~|y_i| \leq 2^{\textstyle 2^{n-2}}$}
\par
\noindent
\centerline{$\forall i \in \{1,\ldots,n\} ~|y_i-1| \leq |x_i-1|$}
\par
\noindent
\centerline{$\forall i,j,k \in \{1,\ldots,n\} ~|y_i+y_j-y_k| \leq |x_i+x_j-x_k|$}
\par
\noindent
\centerline{$\forall i,j,k \in \{1,\ldots,n\} ~|y_i \cdot y_j-y_k| \leq |x_i \cdot x_j-x_k|$}
\vskip 0.2truecm
\par
Since $(\R,+,\cdot,0,1,=,\leq)$ is decidable,
Conjectures 1 and 2 for~$\R$ are decidable for each fixed~$n$.
For a fixed $n$, Conjecture 1 for $\C$ (Conjecture 2 for $\C$) can be translated
into the sentence involving $2n$ real numbers. Since $(\R,+,\cdot,0,1,=,\leq)$ is decidable,
Conjectures 1 and 2 for $\C$ are decidable for each fixed $n$.
\vskip 0.2truecm
\par
Hilbert's tenth problem is to give a computing
algorithm which will tell of a given polynomial equation with integer coefficients
whether or not it has a solution in integers. Yu.~V.~Matijasevi\v{c} proved
(\cite{Matijasevic1}) that there is no such algorithm, see also
\cite{Matijasevic2}, \cite{Davis1973}, \cite{Davis1982}, \cite{Jones1991}.
It implies that Conjecture 1 is false for $\Z$ instead of $\R$~$(\C)$.
Moreover, Matijasevi\v{c}'s theorem implies that Conjecture 1 for $\Z$
is false with any other computable estimation instead of $2^{\textstyle 2^{n-2}}$, so each~$\chi$
in item {\bf (1)} is not computable.
\vskip 0.2truecm
\par
As we have proved Conjecture 1 for $\Z$ is false. We describe a counterexample showing that
Conjecture 1 for $\Z$ is false with $n=21$.
\vskip 0.2truecm
\par
\noindent
{\bf Lemma 2}~([8,~Lemma 2.3,~p.~451]). For each $x \in \Z \cap [2,\infty)$
there exists $y \in \Z \cap [1,\infty)$ such that $1+x^3(2+x)y^2$ is a square.
\newpage
\par
\noindent
{\bf Lemma 3}~([8,~Lemma 2.3,~p.~451]). For each $x \in \Z \cap [2,\infty)$, $y \in \Z \cap [1,\infty)$,
if $1+x^3(2+x)y^2$ is a square, then $y \geq x+x^{x-2}$.
\vskip 0.2truecm
\par
Let us consider the following system over $\Z$. This system consists of two subsystems.
\begin{description}
\item{$(\bullet)$}
\mbox{$~~~x_1=1~~$}
\mbox{$~~x_1+x_1=x_2~~$}
\mbox{$~~x_2 \cdot x_2=x_3~~$}
\mbox{$~~x_3 \cdot x_3=x_4~~$}\\
\mbox{$~~x_4 \cdot x_4=x_5~~$}
\mbox{$~~x_5 \cdot x_5=x_6~~$}
\mbox{$~~x_6 \cdot x_6=x_7~~$}
\mbox{$~~x_6 \cdot x_7=x_8~~$}\\
\mbox{$~~x_2+x_6=x_9~~$}
\mbox{$~~x_8 \cdot x_9=x_{10}~~$}
\mbox{$~~x_{11} \cdot x_{11}=x_{12}~~$}
\mbox{$~~x_{10} \cdot x_{12}=x_{13}~~$}\\
\mbox{$~~x_1+x_{13}=x_{14}~~$}
\mbox{$~~x_{15} \cdot x_{15}=x_{14}~~$}
\item{$(\diamond)$}
\mbox{$~~~x_{16}+x_{16}=x_{17}~~$}
\mbox{$~~x_1+x_{18}=x_{17}~~$}
\mbox{$~~x_{16}+x_{18}=x_{19}~~$}
\mbox{$~~x_{18} \cdot x_{19}=x_{20}~~$}\\
\mbox{$~~x_{12} \cdot x_{21}=x_{20}~~$}
\end{description}
Since $x_1=1$ and $x_{12}=x_{11} \cdot x_{11}$, the subsystem marked with $(\diamond)$ is equivalent to
\par
\noindent
\centerline{$x_{21} \cdot x_{11}^2=(2x_{16}-1)(3x_{16}-1)$}
\par
\noindent
The subsystem marked with $(\bullet)$ is equivalent to 
\par
\noindent
\centerline{$x_{15}^2=1+(2^{16})^3 \cdot (2+2^{16}) \cdot x_{11}^2$}
\par
\noindent
By Lemma 2 the last equation has a solution $(x_{11},x_{15}) \in {\Z}^2$ such that $x_{11} \geq 1$.
By Lemma 1 we can find integers $x_{16}$, $x_{21}$ satisfying $x_{21} \cdot x_{11}^2=(2x_{16}-1)(3x_{16}-1)$.
Thus, the whole system is consistent over $\Z$.
\vskip 0.2truecm
\par
If $(x_1,\ldots,x_{21}) \in {\Z}^{21}$ solves the whole system, then
$x_{15}^2=1+(2^{16})^3 \cdot (2+2^{16}) \cdot |x_{11}|^2$ and
$x_{21} \cdot |x_{11}|^2=(2x_{16}-1)(3x_{16}-1)$.
Since $2x_{16}-1 \neq 0$ and $3x_{16}-1 \neq 0$, $|x_{11}| \geq 1$. By Lemma 3
$$
|x_{11}| \geq 2^{16}+(2^{16})^{\textstyle 2^{16}-2} >
(2^{16})^{\textstyle 2^{16}-2}=2^{\textstyle 2^{20}-32}>2^{\textstyle 2^{21-2}}
$$
{\bf Theorem 8.} If $\Z$ is definable in $\Q$ by an existential formula, then Conjecture 1 fails for $\Q$.
\vskip 0.2truecm
\par
\noindent
{\it Proof.} If $\Z$ is definable in $\Q$ by an existential formula, then $\Z$ is definable in $\Q$
by a Diophantine formula. Let
$$
\forall x_1 \in \Q ~(x_1 \in \Z \Leftrightarrow \exists x_2 \in \Q \ldots \exists x_m \in \Q ~\Phi(x_1,x_2,\ldots,x_m))
$$
where $\Phi(x_1,x_2,\ldots,x_m)$ is a conjunction of the formulae of the form \hbox{$x_i=1$}, \hbox{$x_i+x_j=x_k$,}
\hbox{$x_i \cdot x_j=x_k$,} where $i,j,k \in \{1,\ldots,m\}$. We find an integer $n$ with $2^n \geq m+10$.
Now we are ready to describe a counterexample to Conjecture 1 for $\Q$, this counterexample uses $n+m+11$ variables.
Considering all equations over $\Q$, we can equivalently write down the system
\begin{displaymath}
\left\{
\begin{array}{rcl}
\Phi(x_1,x_2,\ldots,x_m) & &~~~~(1)\\
x_{m+2}^2=1+\left({2^{\textstyle 2^n}}\right)^3 \cdot (2+2^{\textstyle 2^n}) \cdot x_1^2 & &~~~~(2)\\
x_1 \cdot x_{m+1}=1 & &~~~~(3)
\end{array}
\right.
\end{displaymath}
as a conjunction of the formulae of the form $x_i=1$, $x_i+x_j=x_k$, $x_i \cdot x_j=x_k$,
where $i,j,k \in \{1,\ldots,n+m+11\}$. The system is consistent over $\Q$.
Assume that $(x_1,\ldots,x_{n+m+11}) \in {\Q}^{n+m+11}$ solves the system.
Formula~(1) implies that $x_1 \in \Z$. By this and equation~(2), $x_{m+2} \in \Z$.
Equation~(3) implies that $x_1 \neq 0$, so by \hbox{Lemma 3}
$$
|x_1| \geq 2^{\textstyle 2^n}+(2^{\textstyle 2^n})^{\textstyle 2^{\scriptstyle 2^n}-2}>2^{\textstyle 2^{\textstyle n+2^n}-2^{n+1}} \geq 2^{\textstyle 2^{n+2^{\scriptstyle n}-1}} \geq 2^{\textstyle 2^{n+m+11-2}}
$$
\newline
\rightline{$\Box$}
\vskip 0.2truecm
\par
\noindent
{\bf Question.} For which $n \in \{1,2,3,\ldots\}$ there exists a continuous function
\vskip 0.2truecm
\par
\noindent
\centerline{${\R}^n \ni (x_1,\ldots,x_n) \stackrel{{\textstyle f_n}}{{\textstyle \longrightarrow}} (f_{(n,1)}(x_1,\ldots,x_n),\ldots,f_{(n,n)}(x_1,\ldots,x_n)) \in [-2^{\textstyle 2^{n-2}},~2^{\textstyle 2^{n-2}}]^n$}
\vskip 0.2truecm
\par
\noindent
such that
\vskip 0.2truecm
\par
\noindent
\centerline{$\forall (x_1,\ldots,x_n) \in {\R}^n ~\forall i \in \{1,\ldots,n\} ~(x_i=1 \Rightarrow f_{(n,i)}(x_1,\ldots,x_n)=1)$}
\vskip 0.2truecm
\par
\noindent
$\forall (x_1,\ldots,x_n) \in {\R}^n ~\forall i,j,k \in \{1,\ldots,n\} ~(x_i+x_j=x_k \Rightarrow$
\par
\noindent
\rightline{$f_{(n,i)}(x_1,\ldots,x_n)+f_{(n,j)}(x_1,\ldots,x_n)=f_{(n,k)}(x_1,\ldots,x_n))$}
\vskip 0.2truecm
\par
\noindent
$\forall (x_1,\ldots,x_n) \in {\R}^n ~\forall i,j,k \in \{1,\ldots,n\} ~(x_i \cdot x_j=x_k \Rightarrow$
\par
\noindent
\rightline{$f_{(n,i)}(x_1,\ldots,x_n) \cdot f_{(n,j)}(x_1,\ldots,x_n)=f_{(n,k)}(x_1,\ldots,x_n))$}
\vskip 0.2truecm
\par
\noindent
{\bf Theorem 9.} Such functions exist for $n=1$ and $n=2$.
\vskip 0.2truecm
\par
\noindent
{\it Proof.} Case $n=1$. We define $f_{1}: \R \to [0,1]$ by
\begin{displaymath}
f_{1}(x)=\left\{
\begin{array}{ccl}
0 & {\rm if} & x \in (-\infty,0)\\
x & {\rm if} & x \in [0,1]\\
1 & {\rm if} & x \in (1,\infty)
\end{array}
\right.
\end{displaymath}
Case $n=2$. Let $A$ be a closed subset of a metric space $X$ and let $X$
be a locally convex topological linear space.
The Dugundji theorem (\cite{Dugundji}) states that every continuous
map $f:A \to X$ can be extended continuously to all of $X$ in such a way
that the range of the extension lies in the convex hull of $f(A)$.
Applying the Dugundji theorem we find a continuous function
$f_{2}:{\R}^2 \to [-2,2]^2$ with the following properties:
$$
\forall x,y \in [-2,2] ~f_2(x,y)=(x,y)
$$
\begin{displaymath}
\begin{array}{rlllllll}
\forall x \in (-\infty,-2) & (f_{2}(x,1) &=& (-2,1) &\wedge& f_{2}(1,x) &=& (1,-2))\\
\forall x \in (2,\infty) & (f_{2}(x,1) &=& (2,1) &\wedge& f_{2}(1,x) &=& (1,2))\\
\forall x \in (-\infty,-2) & (f_{2}(x,0) &=& (-2,0) &\wedge& f_{2}(0,x) &=& (0,-2))\\
\forall x \in (2,\infty) & (f_{2}(x,0) &=& (2,0) &\wedge& f_{2}(0,x) &=& (0,2))\\
\forall x \in (-\infty,-1) & (f_{2}(x,2x) &=& (-1,-2) &\wedge& f_{2}(2x,x) &=& (-2,-1))\\
\forall x \in (1,2] & (f_{2}(x,2x) &=& (2-x,4-2x) &\wedge& f_{2}(2x,x) &=& (4-2x,2-x))\\
\forall x \in (2,\infty) & (f_{2}(x,2x) &=& (0,0) &\wedge& f_{2}(2x,x) &=& (0,0))\\
\forall x \in (-\infty,-\sqrt{2}) & (f_{2}(x,x^2) &=& (-\sqrt{2},2) &\wedge& f_{2}(x^2,x) &=& (2,-\sqrt{2}))\\
\forall x \in (\sqrt{2},2] & (f_{2}(x,x^2) &=& (\sqrt{4-x^2},4-x^2) &\wedge& f_{2}(x^2,x) &=& (4-x^2,\sqrt{4-x^2}))\\
\forall x \in (2,\infty) & (f_{2}(x,x^2) &=& (0,0) &\wedge& f_{2}(x^2,x)&=&(0,0))
\end{array}
\end{displaymath}
\par
We propose an effective description of a continuous $f_2:{\R}^2 \to [-2,2]^2$.
We define $\sigma: \R \to [-2,2]$ by
\begin{displaymath}
\sigma(x)=\left\{
\begin{array}{ccl}
-2 & {\rm if} & x \in (-\infty,-2)\\
x & {\rm if} & x \in [-2,2]\\
2 & {\rm if} & x \in (2,\infty)
\end{array}
\right.
\end{displaymath}
Let
$$
T=[-2,2]^2 \cup 
$$
$$
\{(x,y) \in {\R}^2:~y=1 \vee x=1 \vee y=0 \vee x=0 \vee y=2x \vee x=2y \vee y=x^2 \vee x=y^2 \}
$$
Let $\rho:{\R}^2 \setminus T \to (0,\infty)$ be defined by
$$
\rho(x,y)=\frac{1}{|x-\sigma(x)|+|y-\sigma(y)|}+
$$
$$
\frac{1}{|y-1|}+\frac{1}{|x-1|}+\frac{1}{|y-0|}+\frac{1}{|x-0|}+
\frac{1}{|y-2x|}+\frac{1}{|x-2y|}+\frac{1}{|y-x^2|}+\frac{1}{|x-y^2|}
$$
and let $g:{\R}^2 \setminus T \to [-2,2]^2$ be defined by
$$
g(x,y)=\frac{1}{\rho(x,y)} \cdot \Biggl(\frac{f_2(\sigma(x),\sigma(y))}{|x-\sigma(x)|+|y-\sigma(y)|}+
$$
$$
\frac{f_2(x,1)}{|y-1|}+\frac{f_2(1,y)}{|x-1|}+\frac{f_2(x,0)}{|y-0|}+\frac{f_2(0,y)}{|x-0|}+
\frac{f_2(x,2x)}{|y-2x|}+\frac{f_2(2y,y)}{|x-2y|}+\frac{f_2(x,x^2)}{|y-x^2|}+\frac{f_2(y^2,y)}{|x-y^2|}\Biggr)
$$
Let $f_{2}|T$ denote $f_{2}$ restricted to $T$.
The function $f_{2}$ has an exact definition on $T$, and $(f_{2}|T) \cup g: {\R}^2 \to [-2,2]^2$ is continuous.
\newline
\rightline{$\Box$}
\vskip 0.2truecm
\par
Let $\K$ be a ring and let $A \subseteq \K$.
We say that a map $f:A \to \K$ is {\sl arithmetic}
if it satisfies the following conditions:
\par
\noindent
if $1 \in A$ then $f(1)=1$,
\par
\noindent
if $a,b \in A$ and $a+b \in A$ then $f(a+b)=f(a)+f(b)$,
\par
\noindent
if $a,b \in A$ and $a \cdot b \in A$ then $f(a \cdot b)=f(a) \cdot f(b)$.
\par
We call an element $r \in \K$ {\sl arithmetically fixed} if there is
a finite set $A \subseteq \K$ (an {\sl arithmetic neighbourhood} of $r$ inside~$\K$)
with $r \in A$ such that each arithmetic map $f:A \to \K$ fixes $r$,
i.e. $f(r)=r$.
If $\K$ is a field, then any $r \in \K$ is arithmetically
fixed if and only if $\{r\}$ is existentially first-order definable in the language
of rings without parameters, see \cite{Tyszka2}.
Articles \cite{Tyszka1}, \cite{Lettl}, \cite{Tyszka2}
dealt with a description of a situation where for an element
in a field there exists an arithmetic neighbourhood.
Article \cite{Tyszka3} describes various types of arithmetic neighbourhoods inside~$\Z$ and~$\Q$.
\vskip 0.2truecm
\par
Let $\widetilde{\K}$ denote the set of all $r \in \K$
that are arithmetically fixed.
Let ${\widetilde{\K}}_n$ ($n=1,2,3,\ldots$) denote the set of
all $r \in \K$ for which there exists an arithmetic neighbourhood
$A$ of $r$ such that ${\rm card}(A) \leq n$. Obviously,
${\widetilde{\K}}_1=\{0,1\}$ and
${\widetilde{\K}}_2 \in \left\{\left\{0,1\right\},\left\{0,1,2\right\},\left\{0,1,2,\frac{1}{2}\right\}\right\}$.
\par
\vskip 0.2truecm
By Theorem~3 in \cite{Tyszka1}~${\widetilde{\R}}_n
\subseteq {\R}^{\rm alg}=\{x \in \R: x {\rm ~is~algebraic~over~}\Q\}$.
By this, Conjecture 1 implies
${\widetilde{\R}}_n \subseteq {\R}^{\rm alg} \cap [-2^{\textstyle 2^{n-2}},~2^{\textstyle 2^{n-2}}]$.
By Corollary~2 in \cite{Tyszka1}~${\widetilde{\C}}_n \subseteq \Q$.
By this, Conjecture 1 implies
${\widetilde{\C}}_n \subseteq \Q \cap [-2^{\textstyle 2^{n-2}},~2^{\textstyle 2^{n-2}}]$.
\vskip 0.2truecm
\par
\noindent
{\bf Theorem 10}~(cf. \cite{Tyszka1}). For each $n \in \{3,4,5,\ldots\}$
we have ${\rm card}({\widetilde{\K}}_n) \leq (n+1)^{n^2+n}+2$.
\vskip 0.2truecm
\par
\noindent
{\it Proof.} If ${\rm card}(\K)<n$ then
${\rm card}(\widetilde{\K}_n) \leq {\rm card}(\K)<n<(n+1)^{n^2+n}+2$.
In the rest of the proof we assume that ${\rm card}(\K) \geq n$.
Let $r \in {\widetilde{\K}}_n \setminus \{0,1\}$ and $A$ is a neighbourhood of $r$
such that ${\rm card}(A) \leq n$. Then each set $B$ with
$A \subseteq B \subseteq \K$ and ${\rm card}(B)=n$ is
a neighbourhood of~$r$.
Observe that $1 \in B$, because in the opposite case the arithmetic map $B \to \{0\}$
moves $r \neq 0$, which is impossible.
Since $r \neq 1$, we can choose $B$ with $\K \supseteq B=\{x_1,\ldots,x_n\} \supseteq A$, where
$x_1=r$, $x_n=1$, and $x_i \neq x_j$ if $i \neq j$. We choose all formulae
$x_i+x_j=x_k$, $x_i \cdot x_j=x_k$
($1 \leq i \leq j \leq n$, $1 \leq k \leq n$) that are satisfied
in~$B$. Joining these formulae with conjunctions we get
some formula~$\Phi$. Let $V$ denote the set of variables
in $\Phi \wedge (x_n=1)$. Observe that $x_1 \in V$, since otherwise for any $s \in \K \setminus \{r\}$
the mapping $f={\rm id}(B \setminus\{r\}) \cup \{(r,s)\}$
satisfies conditions {\bf (1)}-{\bf (3)} and $f(r) \neq r$.
The formula
$\underbrace{\ldots~\exists x_i~\ldots}_{\textstyle {x_i \in V,~i \neq 1}} (\Phi \wedge (x_n=1))$~~~
is satisfied in~$\K$ if and only if $x_1=r$.
\vskip 0.4truecm
\par
For each $(i,j) \in \{(i,j): 1 \leq i \leq j \leq n\}$ there are $n+1$
possibilities:
$$
x_i+x_j=x_1,~~\ldots,~~x_i+x_j=x_n,~~x_i+x_j \not\in \{x_1,\ldots,x_n\}.
$$
\par
For each $(i,j) \in \{(i,j): 1 \leq i \leq j \leq n\}$ there are $n+1$
possibilities:
$$
x_i \cdot x_j=x_1,~~\ldots,~~x_i \cdot x_j=x_n,~~x_i \cdot x_j \not\in \{x_1,\ldots,x_n\}.
$$
Since ${\rm card}(\{(i,j): 1 \leq i \leq j \leq n\})=\frac{n^2+n}{2}$~, the
number of possible formulae $\Phi \wedge (x_n=1)$ does not exceed
$(n+1)^\frac{n^2+n}{2} \cdot (n+1)^\frac{n^2+n}{2}=(n+1)^{n^2+n}$.
Thus \\ ${\rm card}({\widetilde{\K}}_n \setminus \{0,1\}) \leq (n+1)^{n^2+n}$,
so ${\rm card}({\widetilde{\K}}_n) \leq (n+1)^{n^2+n}+2$.
\newline
\rightline{$\Box$}
\vskip 0.2truecm
\par
As we have seen in the proof of Theorem 10,
each $n$-element arithmetic neighbourhood of $r \in \K$
determines a system of equations belonging to
some non-empty subset of~$E_n$. In the ring $\K$,
for each solution of this system the value of variable
$x_1$ is $r$.
\vskip 0.2truecm
\par
Considering all systems $S \subseteq E_i$ $(i=1,2,3)$ we get:
${\widetilde{\Q}}_1={\widetilde{\R}}_1={\widetilde{\C}}_1=\{0,1\}$,
${\widetilde{\Q}}_2={\widetilde{\R}}_2={\widetilde{\C}}_2=\{0,1,2,\frac{1}{2}\}$,
${\widetilde{\Q}}_3={\widetilde{\R}}_3={\widetilde{\C}}_3=
\{0,1,2,\frac{1}{2},-1,3,4,-\frac{1}{2},\frac{1}{4},\frac{3}{2},-2,\frac{1}{3},\frac{2}{3}\}$.
\vskip 0.2truecm
\par
For any ring $\K$ and any $r \in \widetilde{\K}$ we define
$\omega(r) \in \{1,2,3,\ldots\}$ as
\par
\noindent
\centerline{${\rm min} \left\{{\rm card}(A):~ \{r\} \subseteq A \subseteq \K \wedge A {\rm ~is~an~arithmetic~neighbourhood~of~} r {\rm ~inside~} \K \right\}$}
\par
\noindent
As a corollary of Theorem 10 we obtain
$$
\forall n \in \{3,4,5,\ldots\} ~\forall B \subseteq \widetilde{\K} ~\left({\rm card}(B)>(n+1)^{n^2+n}+2 \Longrightarrow \exists r \in B ~\omega(r)>n \right)
$$
\par
Obviously, each ${\widetilde{\K}}_n$ is finite (Theorem 8 gives a concrete
upper bound for ${\rm card}(\widetilde{{\K}}_n)$), so for any subring $\K \subseteq \C$
there exists $\lambda:\{1,2,3,\ldots\} \to \{1,2,3,\ldots\}$ such that
$$
\forall n \in \{1,2,3,\ldots\}
~\forall z \in {\widetilde{\K}}_n ~|z| \leq \lambda(n)
$$
The author does not know whether for $\K=\Z$ there exists a computable
$\lambda:\{1,2,3,\ldots\} \to \{1,2,3,\ldots\}$ with the above property.
\vskip 0.2truecm
\par
\noindent
{\bf Conjecture 3.} Let $\G$ be an additive subgroup of $\C$. Let
$S$ be a consistent system of equations in $x_1,x_2,\ldots,x_n \in \G$, where each equation in $S$ is one
of the following two forms: $x_i=1$ or $x_i+x_j=x_k$. Then $S$ has
a solution $(x_1,x_2,\ldots,x_n) \in {(\G \cap \Q)}^n$ in which $|x_j| \leq 2^{n-1}$ for each $j$.
\vskip 0.2truecm
\par
In case when $\G \supseteq \Q$ we will prove a weaker version of Conjecture 3 with the estimation given by $(\sqrt{5})^{n-1}$.
\vskip 0.2truecm
\par
\noindent
{\bf Observation 3.} If ${\cal A} \subseteq {\C}^k$ is an affine
subspace and ${\rm card}~{\cal A}>1$, then there exists $m \in \{1,2,\ldots,k\}$ with
$$
\emptyset \neq {\cal A} \cap \{(x_1,x_2,\ldots,x_k) \in {\C}^k:~x_m+x_m=x_m\} \subsetneq {\cal A}
$$
{\bf Theorem 11.} Let $S$ be a consistent system of equations in complex
numbers $x_1,x_2,\ldots,x_n$, where each equation in $S$ is one
of the following two forms: $x_i=1$ or $x_i+x_j=x_k$. Then $S$ has a rational
solution $(x_1,x_2,\ldots,x_n)$ in which $|x_j| \leq (\sqrt{5})^{n-1}$ for each $j$.
\vskip 0.2truecm
\par
\noindent
{\it Proof.} We shall describe how to find a solution $(x_1,x_2,\ldots,x_n) \in {\Q}^n$ in which $|x_j| \leq (\sqrt{5})^{n-1}$
for each $j$. We can assume that for a certain $i \in \{1,2,\ldots,n\}$ the equation $x_i=1$ belongs to $S$,
as otherwise $(0,0,\ldots,0)$ is a solution.
Without lost of generality we can assume that the equation $x_1=1$ belongs to $S$.
Each equation belonging to $S$ has a form
$$
a_1x_1+a_2x_2+\ldots+a_nx_n=b,
$$
where $a_1,a_2,\ldots,a_n,b \in \Z$. Since $x_1=1$, we can equivalently write this equation as
$$
a_2x_2+a_3x_3+\ldots+a_nx_n=b-a_1
$$
We receive a system of equations whose set of solutions is a non-empty affine subspace
${\cal A} \subseteq {\C}^{n-1}$.
If ${\rm card~} {\cal A}>1$, then by Observation 3 we find $m \in \{2,3,\ldots,n\}$ for which
$$
\emptyset \neq {\cal A} \cap \{(x_2,x_3,\ldots,x_n) \in {\C}^{n-1}:~x_m+x_m=x_m\} \subsetneq {\cal A}
$$
The procedure described in the last sentence is applied to the affine subspace
$$
{\cal A} \cap \{(x_2,x_3,\ldots,x_n) \in {\C}^{n-1}:~x_m+x_m=x_m\}
$$
and repeated until one point is achieved.
The maximum number of procedure executions is $n-1$. The received one-point affine 
subspace is described by equations belonging to a certain set
$$
{\cal U} \subseteq \{x_i=1:~i \in \{2,3,\ldots,n\}\} \cup \{x_i+x_j=x_k:~i,j,k \in \{1,2,\ldots,n\}, i+j+k>3\}
$$
Each equation belonging to ${\cal U}$ has a form
$$
a_2x_2+a_3x_3+\ldots+a_nx_n=c,
$$
where $a_2,a_3,\ldots,a_n,c \in \Z$.
Among these equations, we choose $n-1$ linearly independent equations.
Let ${\bf A}$ be the matrix of the system, and the system of equations has the following form
$$
{\bf A} \cdot
\left[\begin{array}{c}x_2 \\x_3 \\ \vdots \\ x_n
\end{array}\right]
=
\left[\begin{array}{c}c_2 \\c_3 \\ \vdots \\ c_n
\end{array}\right]
$$
Let ${\bf A}_j$ be the matrix formed by replacing the $j$-th column of ${\bf A}$ by the column vector $c_2,c_3,\ldots,c_n$.
Obviously, $\det({\bf A}) \in \Z$, and $\det({\bf A}_j) \in \Z$ for each $j \in \{1,2,\ldots,n-1\}$.
By Cramer's rule
$x_j=\frac{\det({\bf A}_{j-1})}{\det({\bf A})} \in \Q$ for each $j \in \{2,3,\ldots,n\}$.
\vskip 0.15truecm
\par
\noindent
When the row of matrix ${\bf A}$ corresponds to the equation
$x_i=1$~($i> 1$),
then the entries in the row are $1$, $0$ ($n-2$ times),
while the right side of the equation is $1$.
\vskip 0.15truecm
\par
\noindent
When the row of matrix ${\bf A}$ corresponds to the equation
$x_1+x_1=x_i$~($i> 1$),
then the entries in the row are $1$, $0$ ($n-2$ times),
while the right side of the equation is $2$.
\vskip 0.15truecm
\par
\noindent
When the row of matrix ${\bf A}$ corresponds to one of the equations:
$x_1+x_i=x_1$ or $x_i+x_1=x_1$~($i>1$),
then the entries in the row are $1$, $0$ ($n-2$ times),
while the right side of the equation is $0$.
\vskip 0.15truecm
\par
\noindent
When the row of matrix ${\bf A}$ corresponds to one of the equations:
$x_1+x_i=x_j$ or $x_i+x_1=x_j$ ($i>1$, $j>1$, $i \neq j$),
then the entries in the row are $1$, $-1$, $0$ ($n-3$ times),
while the right side of the equation is $1$. 
\vskip 0.15truecm
\par
\noindent
When the row of matrix ${\bf A}$ corresponds to the equation
$x_i+x_i=x_1$~($i>1$),
then the entries in the row are $2$, $0$ ($n-2$ times),
while the right side of the equation is $1$.
\vskip 0.15truecm
\par
\noindent
When the row of matrix ${\bf A}$ corresponds to the equation
$x_i+x_j=x_1$ ($i>1$, $j>1$, $i \neq j$),
then the entries in the row are $1$, $1$, $0$ ($n-3$ times),
while the right side of the equation is $1$.
\vskip 0.15truecm
\par
From now on we assume that $i,j,k \in \{2,3,\ldots,n\}$.
\vskip 0.15truecm
\par
\noindent
When the row of matrix ${\bf A}$ corresponds to the equation
$x_i+x_j=x_k$~($i \neq j$, $i \neq k$, $j \neq k$),
then the entries in the row are $1$, $1$, $-1$, $0$ ($n-4$ times),
while the right side of the equation is $0$.
\vskip 0.15truecm
\par
\noindent
When the row of matrix ${\bf A}$ corresponds to the equation
$x_i+x_i=x_k$~($i \neq k$),
then the entries in the row are $2$, $-1$, $0$ ($n-3$ times),
while the right side of the equation is $0$.
\newpage
\par
\noindent
When the row of matrix ${\bf A}$ corresponds to the equation
$x_i+x_j=x_k$~($k=i$ or $k=j$),
then the entries in the row are $1$, $0$ ($n-2$ times),
while the right side of the equation is $0$.
\vskip 0.2truecm
\par
Contradictory equations, e.g. $x_1+x_i=x_i$ do not belong to ${\cal U}$,
therefore their description has been disregarded.
The presented description shows that each row of matrix ${\bf A}_j$ ($j \in \{1,2,\ldots,n-1\}$)
has the length less than or equal to $\sqrt{5}$.
By Hadamard's inequality $|\det({\bf A}_j)| \leq (\sqrt{5})^{n-1}$ for each $j \in \{1,2,\ldots,n-1\}$.
Hence, $|x_j|=\frac{|\det({\bf A}_{j-1})|}{|\det({\bf A})|} \leq |\det({\bf A}_{j-1})| \leq (\sqrt{5})^{n-1}$
for each $j \in \{2,3,\ldots,n\}$.
\newline
\rightline{$\Box$}
\vskip 0.2truecm
\par
Concerning the above proof, without lost of generality we can assume that all equations $x_i=1$ ($i>1$)
do not belong to $S$. Indeed, if $i>1$ and the equation $x_i=1$ belongs to $S$, then we replace
$x_i$ by $x_1$ in all equations belonging to $S$. In this way the problem reduces to the same problem
with a smaller number of variables. Therefore, for proving Theorem 11 (or any other bound)
it is sufficient to consider only these \hbox{systems $S$} of $n$ equations which have a unique solution $(x_1,\ldots,x_n)$
and contain the equation $x_1=1$ and $n-1$ equations of the form $x_i+x_j=x_k$ ($i,j,k \in \{1,2,\ldots,n\}$).
Let ${\bf B}$ be the matrix of the linear system consisting of the aforementioned $n-1$ equations of the form $x_i+x_j=x_k$.
Let ${\bf A}$ be the matrix of the following linear system
\par
\noindent
\centerline{~~~~~$x_1=1$}
$$
{\bf B} \cdot
\left[\begin{array}{c}x_1 \\x_2 \\ \vdots \\ x_n
\end{array}\right]
=
\left[\begin{array}{c}0 \\0 \\ \vdots \\ 0
\end{array}\right]
$$
and let ${\bf B}_j$ be the matrix formed by deleting the $j$-th column of ${\bf B}$. By Cramer's rule
$|x_j|=|\frac{\det({\bf B}_j)}{\det({\bf A})}| \leq |\det({\bf B}_j)|$ for each $j \in \{2,3,\ldots,n\}$.
By this, the following Conjecture 4 implies Conjecture 3 for the case when $\G \supseteq \Q$.
\vskip 0.2truecm
\par
\noindent
{\bf Conjecture 4.} Let ${\bf B}$ be a matrix with $n-1$ rows and $n$ columns, $n \geq 2$.
Assume that each row of ${\bf B}$, after deleting all zeros, forms a sequence belonging to
\par
\noindent
\centerline{$\{
\langle 1 \rangle,
\langle -1,2 \rangle,
\langle 2,-1 \rangle,
\langle -1,1,1 \rangle,
\langle 1,-1,1 \rangle,
\langle 1,1,-1 \rangle\}$}
\par
\noindent
We conjecture that after deleting any column of ${\bf B}$ we get the matrix
whose determinant has absolute value less than or equal to $2^{n-1}$.
\vskip 0.4truecm
\par
In case when $\G=\R$ ($\G=\C$) Conjecture 5 strengthens Conjecture 3.
\vskip 0.2truecm
\par
\noindent
{\bf Conjecture 5.} For each $x_1,\ldots,x_n \in \R$~($\C$) there exist $y_1,\ldots,y_n \in \R$~($\C$) such that
\par
\noindent
\centerline{$\forall i \in \{1,\ldots,n\} ~|y_i| \leq 2^{n-1}$}
\par
\noindent
\centerline{$\forall i \in \{1,\ldots,n\} ~|y_i-1| \leq |x_i-1|$}
\par
\noindent
\centerline{$\forall i,j,k \in \{1,\ldots,n\} ~|y_i+y_j-y_k| \leq |x_i+x_j-x_k|$}
\vskip 0.2truecm
\par
For a positive integer $n$ we define the set of equations $W_n$ by
\par
\noindent
\centerline{$W_n=\{x_i=1:~1 \leq i \leq n\}~\cup \{x_i+x_j=x_k:~1 \leq i \leq j \leq n,~1 \leq k \leq n\}$}
\par
\noindent
If a system $S \subseteq W_n$ is consistent over $\R$~($\C$) and maximal with respect to inclusion,
then (cf. the proof of Theorem 11) $S$ has a unique rational solution $(\widehat{x_1},\ldots,\widehat{x_n})$ given by Cramer's rule.
Hence,
\par
\noindent
\centerline{Conjecture 3 for $\R$ $\Longleftrightarrow$ Conjecture 3 for $\Q$ $\Longleftrightarrow$ Conjecture 3 for $\C$}
\vskip 0.2truecm
\par
Conjecture 3 holds true for each $n \in \{1,2,3,4\}$ and each additive subgroup $\G \subseteq \C$.
It follows from the following Observation 4.
\vskip 0.2truecm
\par
\noindent
{\bf Observation 4.} Let $n \in \{1,2,3,4\}$, and let $S \subseteq W_n$ be a system that is consistent over
the additive subgroup $\G \subseteq \C$.
If $(x_1,\ldots,x_n) \in {\G}^n$ solves $S$, then $(\widehat{x_1},\ldots,\widehat{x_n})$ solves $S$,
where each $\widehat{x_i}$ is suitably chosen from $\{x_i,0,1,2,\frac{1}{2}\} \cap \{z \in \G:~|z| \leq 2^{n-1}\}$.
\vskip 0.2truecm
\par
Multiple execution of the algorithm described in items
{\bf (22)}-{\bf (23)} yields partial (as probabilistic)
resolution of Conjecture 3 for $\R$ and $n \geq 2$.
This algorithm resolves Conjecture 3 for some randomly chosen subsystem of $W_n$.
\vskip 0.2truecm
\par
\noindent
{\bf (22)}~~We define by induction a finite sequence $(s_1,\ldots,s_n)$
of equations belonging to $W_n \setminus \{x_i=1:~1 \leq i \leq n\}$. As $s_1$ we put the equation $x_1+x_1=x_1$.
When $1 \leq i \leq n-1$ and the sequence $(s_1,\ldots,s_i)$ is defined,
then as $s_{i+1}$ we put the first randomly chosen $h \in W_n \setminus \{x_i=1:~1 \leq i \leq n\}$
for which the equations $s_1,\ldots,s_i,h$ are linearly independent.
\vskip 0.2truecm
\par
\noindent
{\bf (23)}~~We resolve Conjecture 3 for $\R$ for the system $\{x_1=1,s_2,\ldots,s_n\}$.
This system has a unique solution $(1,\widehat{x_2},\ldots,\widehat{x_n})$ given by Cramer's rule,
$\widehat{x_2},\ldots,\widehat{x_n} \in \Q$.
\vskip 0.2truecm
\par
The algorithm terminates with probability $1$.
For $n=5$, the following code in {\sl MuPAD} performs the algorithm with $1000$ iterations.
\begin{quote}
\begin{verbatim}
SEED:=time():
r:=random(1..5):
idmatrix:=matrix::identity(5):
u:=linalg::row(idmatrix,i) $i=1..5:
max_norm:=1:
for k from 1 to 1000 do
a:=linalg::row(idmatrix,1):
rank:=1:
while rank<5 do
m:=matrix(u[r()])+matrix(u[r()])-matrix(u[r()]):
a1:=linalg::stackMatrix(a,m):
rank1:=linalg::rank(a1):
if rank1 > rank then a:=linalg::stackMatrix(a,m) end_if:
rank:=linalg::rank(a):
end_while:
x:=(a^-1)*linalg::col(idmatrix,1):
max_norm:=max(max_norm,norm(x)):
print(max_norm):
end_for:
\end{verbatim}
\end{quote}
For another algorithm, implemented in {\sl Mathematica}, see \cite{Kozlowski}.
\vskip 0.2truecm
\par
In case when $\G=\Z$ we will prove a weaker version of Conjecture 3 with the estimation given by $(\sqrt{5})^{n-1}$.
\vskip 0.2truecm
\par
\noindent
{\bf Lemma 4} (\cite{Borosh}). Let ${\bf A}$ be a matrix with $m$ rows, $n$ columns, and integer entries.
Let $b_1,b_2,\ldots,b_m \in \Z$, and the matrix equation
$$
{\bf A} \cdot
\left[\begin{array}{c}x_1 \\x_2 \\ \vdots \\ x_n
\end{array}\right]
=
\left[\begin{array}{c}b_1 \\b_2 \\ \vdots \\ b_m
\end{array}\right]
$$
defines the system of linear equations with rank $m$. Denote by $\delta$ the maximum
of the absolute values of the $m \times m$ minors of the augmented
matrix $({\bf A},b)$. We claim that if the system is consistent over $\Z$,
then it has a solution in $({\Z } \cap [-\delta,\delta])^n$.
\vskip 0.2truecm
\par
\noindent
{\bf Theorem 12.} Let $S$ be a consistent system of equations in integers
$x_1,x_2,\ldots,x_n$, where each equation in $S$ is one
of the following two forms: $x_i=1$ or $x_i+x_j=x_k$. Then $S$ has an integer
solution $(x_1,x_2,\ldots,x_n)$ in which $|x_j| \leq (\sqrt{5})^{n-1}$ for each $j$.
\vskip 0.2truecm
\par
\noindent
{\it Proof.} We shall describe how to find a solution $(x_1,x_2,\ldots,x_n) \in {\Z}^n$ in which $|x_j| \leq (\sqrt{5})^{n-1}$
for each $j$. We can assume that for a certain $i \in \{1,2,\ldots,n\}$ the equation $x_i=1$ belongs to $S$,
as otherwise $(0,0,\ldots,0)$ is a solution.
Without lost of generality we can assume that the equation $x_1=1$ belongs to $S$.
Analogously as in the proof of Theorem 11, we construct a system of linear equations
with variables $x_2,\ldots,x_n$. For the augmented matrix of this system, the length of each row
is not greater than $\sqrt{5}$. We finish the proof by applying Hadamard's inequality and Lemma 4.
\newline
\rightline{$\Box$}

Apoloniusz Tyszka\\
Technical Faculty\\
Hugo Ko{\l}{\l}\c{a}taj University\\
Balicka 116B, 30-149 Krak\'ow, Poland\\
E-mail: {\it rttyszka@cyf-kr.edu.pl}
\end{document}